\documentclass[12pt, reqno]{amsart}
\usepackage{amsmath}
\usepackage{dsfont}
\usepackage{amsfonts}
\usepackage{amssymb}
\usepackage{aliascnt}
\usepackage{graphicx}
\usepackage{mathrsfs}
\usepackage{enumerate}
\usepackage[bookmarks=true,pdfstartview=FitH, pdfborder={0 0 0},
colorlinks=true,citecolor=blue, linkcolor=blue]{hyperref}

\newtheorem{theorem}{Theorem}[section]
\newaliascnt{conj}{theorem}
\newaliascnt{cor}{theorem}
\newaliascnt{lemma}{theorem}
\newaliascnt{fact}{theorem}
\newaliascnt{claim}{theorem}
\newaliascnt{prop}{theorem}
\newaliascnt{definition}{theorem}

\newtheorem{prop}[prop]{Proposition}

\aliascntresetthe{conj} \aliascntresetthe{cor}
\aliascntresetthe{lemma} \aliascntresetthe{fact}
\aliascntresetthe{claim} \aliascntresetthe{prop}
\aliascntresetthe{definition}

\theoremstyle{definition}
\newaliascnt{example}{theorem}

\aliascntresetthe{example}

\theoremstyle{remark}
\newaliascnt{rmk}{theorem}

\aliascntresetthe{rmk} 

\def\sek~{\S{}}

\numberwithin{equation}{section}
\newcommand{\curv}{\operatorname{curv}}

\newcommand{\diam}{\operatorname{diam}}

\newcommand{\dis}{\operatorname{d}}

\newcommand{\dbl}{\operatorname{dbl}}

\newcommand{\RR}{\mathds{R}}
\newcommand{\ZZ}{\mathds{Z}}

\renewcommand{\SS}{\mathbf{S}}
\newcommand{\TT}{\mathbf{T}}
\newcommand{\KK}{\mathbf{K}}
\newcommand{\RP}{\mathds{RP}}

\begin{document}
%\linenumbers
\title[4-dimensional Soul Theorem]
{Soul Theorem for 4-dimensional Topologically Regular Open
Nonnegatively Curved Alexandrov Spaces}
\author{Jian GE}
\thanks{This paper has been accepted for publication in
Proc AMS}

\email{jge@nd.edu}

\address{Mathematics Department \\University of Notre Dame\\ Notre Dame, IN 46556, USA}
\maketitle
\begin{abstract}
In this paper, we study the topology of topologically regular
4-dimensional open non-negatively curved Alexandrov spaces. These
spaces occur naturally as the blow-up limits of compact Riemannian
manifolds with lower curvature bound. These manifolds have also been
studied by Yamaguchi in his preprint \cite{Yam2002}. Our main tools
are gradient flows of semi-concave functions and critical point
theory for distance functions, which have been used to study the
$3$-dimensional collapsing theory in the paper \cite{CaoG2010}. The
results of this paper will be used in our future studies of
collapsing 4-manifolds, which will be discussed elsewhere.
\end{abstract}
\section{Introduction}\label{sec:0:Intro}
The topology of noncompact manifold with a complete metric of
nonnegative sectional curvatures was studied by Gromoll-Meyer in
\cite{GM1969} and Cheeger-Gromoll in \cite{CG1972}.

\begin{theorem}[Soul Theorem, \cite{CG1972}]\label{thm:intr:CGsoul}
Let $M^n$ be an $n$-dimensional noncompact manifold with a complete
metric of nonnegative sectional curvature. Then there exists a
compact totally geodesic embedded submanifold $S\subset M$ with
nonnegative sectional curvature, such that $M^n$ is diffeomorphic to
the normal bundle $\nu(S)$ of $S$ in $M^n$.
\end{theorem}

If in addition to the assumption above, there exists $p\in M^n$ such
that the sectional curvatures at $p$ are all positive, then $M^n$ is
diffeomorphic to $\RR^n$. This is called Cheeger-Gromoll soul
conjecture. It was proved by Perelman in \cite{Per1994}.

For Alexandrov spaces, (see \cite{BGP1992} and \cite{BBI2001} for
basics of Alexandrov spaces) Perelman proved a similar result:
\begin{theorem}[Soul Theorem for Alexandrov space, \cite{Per1991}]\label{thm:intr:PerSoul}
Let $X^n$ be an $n$-dimensional non-compact Alexandrov space with
nonnegative curvature, then there exits a closed totally convex
subset $S\subset X^n$, such that $S$ is a deformation retraction of
$X^n$.
\end{theorem}

Similar to the soul conjecture, Cao-Dai-Mei \cite{CDM2007,CDM2009}
proved that if in addition to the conditions of
\autoref{thm:intr:PerSoul}, one assumes that $X^n$ has positive
curvature in a metric ball, then $X^n$ is contractible. Unlike the
manifold case, \autoref{thm:intr:PerSoul} is the best topological
result one can expect, in the sense that in general $X^n$ is not
homeomorphic to the normal bundle of $S$.

By Perelman's stability theorem, if Alexandrov space is the limit of
sequence of Riemannian manifolds with lower curvature bound, then
it's a topological manifold. In fact Kapovitch showed in
\cite{Kap2002}
\begin{theorem}[\cite{Kap2002}]\label{thm:intr:KapSmooth}
If $X^n$ is the limit of a sequence of $n$-dimensional Riemannian
manifold with the same lower curvature bound $k$, then $\Sigma_pX^n$
is homeomorphic to $(n-1)$-sphere $\SS^{n-1}$ for any $p\in X^n$.
Moreover all the iterated space of directions are homeomorphic to
spheres.
\end{theorem}

It's still unknown whether an Alexandrov space, which satisfies the
conclusion of \autoref{thm:intr:KapSmooth}, can be realized as a
limit of Riemannian manifolds with the same dimension and same lower
curvature bound. In this paper, we consider the class of
4-dimensional topologically regular open nonnegatively curved
Alexandrov spaces, in the sense that all space of directions are
spheres. These Alexandrov spaces occurs naturally as the blow-up
limits of compact Riemannian manifolds with lower curvature bound,
thus play an important role in the study of collapsing under a lower
curvature bound. We will prove the following
\begin{theorem}[Main Theorem]\label{thm:intr:Main}
Let $X^4$ be as above, $S$ be a soul of $X^4$, then $X^4$ is
homeomorphic a open disk bundle over $S$:
\begin{enumerate}[{\rm (1)}]
\item
If $\dim S=0$, then $X^4$ is homeomorphic to $\RR^4$;
\item
If $\dim S=1$, then $X^4$ is homeomorphic to $\RR^3$ bundle over
$\SS^1$;
\item
If $\dim S=2$, then $X^4$ is homeomorphic to $\RR^2$ bundle over
$S$, where $S=\SS^2, \RP^2, \TT^2$ or $\KK^2$.
\item
If $\dim S=3$, then $X^4$ is homeomorphic to line bundle over $S$,
where $S=\SS^3/\Gamma, \TT^3/\Gamma, (\SS^2\times\SS^1)/\Gamma$, and
$\Gamma$ some subgroup of isometric group of $S$ acting freely on
$S$.
\end{enumerate}
\end{theorem}

This theorem will be used to study the collapsing 4-manifolds, which
will discussed elsewhere. \autoref{thm:intr:Main} has also been
studied in the preprint \cite{Yam2002} Chap 15,16 in a traditional
way. Our main tools are the gradient flow of semi-concave functions
and Perelman's version of Fibration theorem, which have been used
extensively in \cite{CaoG2010} to study the 3-dimensional collapsing
manifolds under a lower curvature bound.

\section{Construction of Soul}
In this section we recall Cheeger-Gromoll's construction of soul for
$X^n$, where $X^n$ is an $n$-dimensional non-negatively curved open
complete Alexandrov space. (c.f. \cite{CG1972}, \cite{Per1991}).

Fix $p\in X^n$, the Busemann function can be defined by
\begin{equation}
b(x)=\lim_{t\to \infty}[\dis (x, \partial(B(p, t)))-t].
\end{equation}
where $B(p,t)$ is the ball centered at $p$ with radius $t$. Denote
the super-level set $b^{-1}([a, +\infty))$ by $\Omega^a$. We have
\begin{prop}[\cite{CG1972}, \cite{Per1991}, \cite{Wu1979}]\label{prop:con:Filtration}
Let $X^n$, $\Omega^a$ be as above, then the following hold
\begin{enumerate}[{\rm (1)}]
\item The Busemann function $b$ is concave and bounded above.

\item
$\Omega^a$ is compact and totally convex for all $a\le
a_0:=\max_{x\in X^n} b(x)$.

\item
$a<b\le a_0$ implies $\Omega^b\subset \Omega^a$ and
$$\Omega^b=\{x\in \Omega^a| \dis(x, \partial \Omega^a)\le b-a\}$$

\item
There is a filtration of $\Omega(0):=\Omega^{a_0}$ by totally convex
sets:
$$\Omega(0)\supset \Omega(1)\supset \cdots \supset\Omega(k)$$
where $\Omega(i+1)$ is the maximum set of $f_i(x)=\dis_{\Omega(i)}
(x,\partial \Omega(i))$, and $\partial\Omega(k)=\varnothing$.
\end{enumerate}
\end{prop}

We call $S=\Omega(k)$ a soul of the type $(s, m)$, if the dimension
of the soul is $s$ and the dimension of $\Omega(0)$ is $m$.\\

We call a geodesic $\gamma: (-\infty, +\infty)\to X^n$ a line in a
metric space $X^n$, if $\dis(\gamma(t), \gamma(s))= |t-s|$ for
$\forall t, s\in \RR$. The splitting theorem reduces our discussion
of $4$-dimension to 3-dimension when $X^4$ admits a line.
\begin{theorem}[\cite{GP1989}]
Let $X^n$ be open non-negatively curved Alexandrov space and assume
that $X^n$ admits a line, then $X^n$ splits isometrically as $X^n=
N^{n-1} \times \RR$, where $N^{n-1}$ is a $(n-1)$ dimensional open
non-negatively curved Alexandrov space.
\end{theorem}

In order to handle the non-smooth metric Alexandrov space,
Perelman's Stability Theorem and his version of Fibration Theorem
are extensively used in this paper. Let's recall these results.

\begin{theorem}[Stability Theorem \cite{Per1991}, \cite{Kap2007}]\label{thm:Stability}
Let $\{X^n_{\alpha}\}_{\alpha=1}^{\infty}$ be a sequence of
$n$-dimensional Alexandrov spaces with $\curv\ge -1$ converging to
an Alexandrov space with same dimension: $\lim_{\alpha\to \infty}
X^n_{\alpha}=X^n$. Then $X^n_{\alpha}$ is homeomorphic to $X^n$ for
$\alpha$ large.
\end{theorem}

The stability theorem for pointed spaces can be stated in a
similarly way. The fibration theorem states that

\begin{theorem}[Fibration Theorem \cite{Per1991, Per1993}]
Let $X^n$ be an $n$-dimensional Alexandrov space, $U$ a domain in
$X^n$, $f: U\to \RR^k$ be an admissible function, having no critical
point and proper on $U$, then it's restriction to this domain is a
locally trivial fiber bundle.
\end{theorem}

We refer to \cite{Per1991} for the definitions of admissible
functions and regular map.

\section{Soul theorem for 3-dimensional Alexandrov space}
The topology of 3-dimensional open non-negatively curved Alexandrov
space was studied in \cite{SY2000}, (cf also \cite{CaoG2010}).

\begin{theorem}[\cite{SY2000}]\label{thm:3D:3DSoul}
Let $X^3$ be an open complete $3$-manifold with a possibly singular
metric of non-negative curvature. Suppose that $X^3$ is oriented and
$S^s$ is a soul of $X^3$. Then the following is true.
\begin{enumerate}[{\rm (1)}]
\item
When $\dim(S^s) = 1$, then the soul of $X^3$ is isometric to a
circle. Moreover, its universal cover $\tilde X^3$ is isometric to
$\tilde X^2\times \RR$, where $\tilde X^2$ is homeomorphic to $\Bbb
R^2$;
\item
When $\dim(S^s) = 2$, then the soul of $X^3$ is homeomorphic to
$\SS^2/\Gamma$ or $\TT^2/\Gamma$. Furthermore, $X^3$ is isometric to
one of four spaces: $\SS^2 \times \RR$, $\RP^2 \ltimes \RR = (\SS^2
\times \RR)/ \ZZ_2$, $\TT^2 \times \RR$ or $\KK^2 \ltimes \RR =
(\TT^2 \times \RR)/\ZZ_2$ and $\RP^2 \ltimes \RR$ which is
homeomorphic to $ [\RP^3 - D^3]$;

\item
When $\dim(S^s) = 0$, then the soul of $X^3$ is a single point and
$X^3$ is homeomorphic to $\RR^3$.

\end{enumerate}
\end{theorem}

Throughout this paper $\SS^n$ denotes the standard $n$-sphere,
$\TT^n$ denotes the $n$-dimensional torus, $\KK^2$ denotes the Klein
bottle and $D^3$ denotes the standard $3$-ball.

\section{Proof of the main theorem}
A key observation by K. Grove is that the distance function to the
soul has no critical point in $X^n-S^s$ when $X^n$ is a smooth
Riemannian manifold. However this is no longer true for Alexandrov
spaces even for topologically regular one. For example let $L$ be
the closed half strip $\{(x,y)\in \RR^2|x\ge 0, 0\le y\le 1\}$, then
the double of $L$ which denoted by $\dbl (L)$ is a 2-dimensional
Alexandrov space with non-negative curvature, which is homeomorphic
to $\RR^2$ and has point soul $(0, \frac12)$, it's clear that the
distance function has critical points $(0,0)$ and $(0, 1)$. However
for Alexandrov space, \cite{CaoG2010} derived a modified result
similar to Grove's observation.

\begin{prop}[\cite{CaoG2010} Proposition 2.5]\label{lem:4D:nonCritcal}
The function $f(x)= \dis_{X^n}(x, A)$ has no critical point for
$x\in [X^n-\Omega(0)]$, where $A\subset \Omega(0)$.
\end{prop}

\begin{figure*}[ht]
% Requires \usepackage{graphicx}
\includegraphics[width=150pt]{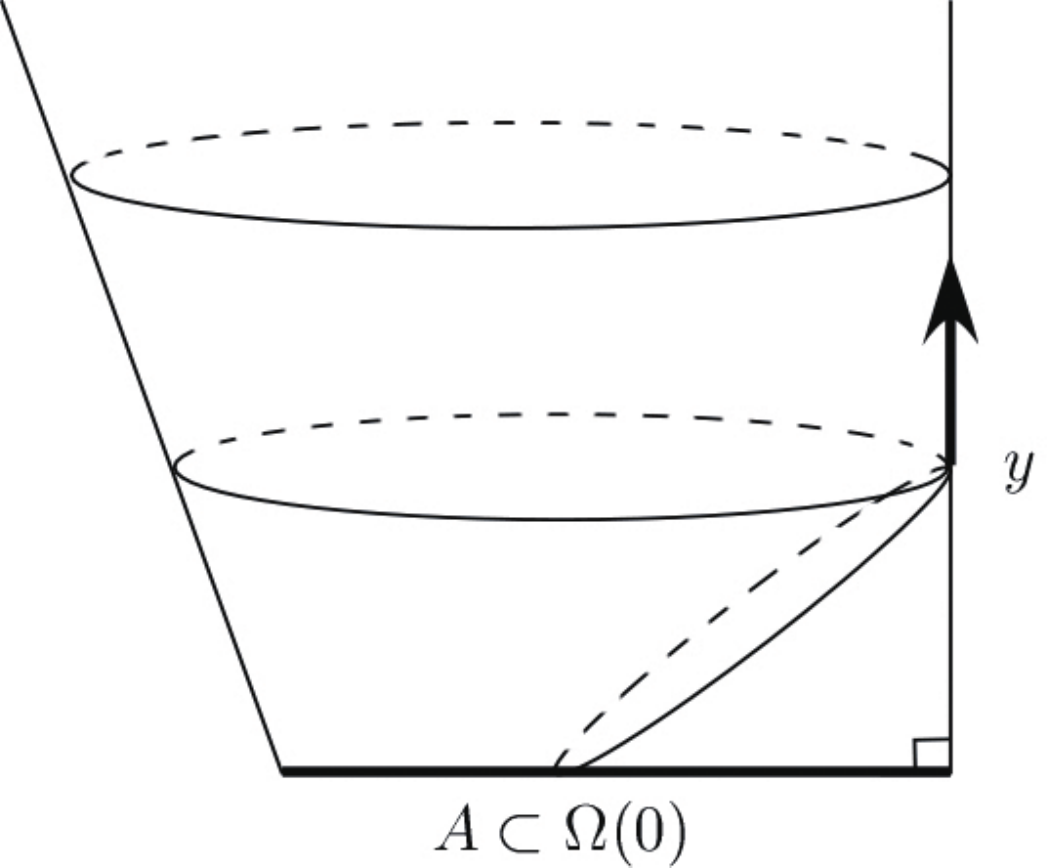}\\
\caption{$\dis_{X^n}(x, A)$ has no critical point in
$[X^n-\Omega(0)]$}\label{Pic:0}
\end{figure*}

For simplicity we assume that $a_0=0$ by adding a constant if
needed. Let $A=\Omega(0)$. Using \autoref{lem:4D:nonCritcal}, we see
that $f(x)=\dis(x, A)$ has no critical value within $(\varepsilon,
+\infty)$. It follows by Perelman's Fibration Theorem, that
$X^n\cong B_{X^n}(\Omega(0), \varepsilon)$ for $\varepsilon>0$,
where $B_X(A,r)$ denoted the set of points with distance $\le r$ to
set $A$ in metric $X$.

We will divide the proof of our Main Theorem into the flowing cases.

\subsection{Soul of the type $(s=3, m=3)$,
$X^4 = S\times \RR^1$ or $(S\times \RR^1)/\ZZ_2$ }\label{sub:4D:1}
\begin{proof}(cf. the proof of Theorem 2.21 in \cite{CaoG2010})
In this case, $S=\Omega(0)$ and has dimension $3$. Since $X^4$ is
topologically regular, $S$ is a topological manifold and
non-negatively curved, hence by Hamilton's classification of
3-dimensional manifolds with non-negatively Ricci curvature, $S$ is
homeomorphic to $\SS^3/\Gamma$, $\TT^3/\Gamma$ or $\SS^2\times\SS^1/
\Gamma$, where $\Gamma$ is a subgroup of the isometry group of
$\SS^3, \TT^3, \SS^2\times \SS^1$.

For $p\in S$, we know $\Sigma_p(S)$ is homeomorphic to $\SS^2$,
which divide $\Sigma_p(X^4)\cong \SS^3$ into two parts, denoted by
$A^{\pm}$. Since $S$ is totally convex, $\Sigma_p(S)$ is convex in
$\Sigma_p(X^4)$, therefore $r_{\Sigma_p(S)}|_{A^{\pm}}$ have a
unique maximum point $\xi^{\pm}$ in $A^{\pm}$. Denote the maximum
values by $\ell^{\pm}$, i.e. $r_{\Sigma_p(S)}(\xi^{\pm})=
\ell^{\pm}$. Since $\Omega(0)=S$ is the set of maximum points for
Busemann function, by the first variation theorem, we know
$\ell^{\pm}\ge \pi/2$. On the other hand if $\gamma: [0, \ell]\to
\Sigma_p(X^4)$ is a shortest geodesic connecting $\xi^-$ to $\xi^+$,
and let $t_0\in [0, \ell]$ satisfying $\gamma(t_0)\in \Sigma_p(S)$.
By triangle inequality we know

\begin{equation} \label{eq:1}
\begin{split}
\dis(\xi^-, \xi^+)&= \dis(\xi^-, \gamma(t_0))+ \dis(\gamma(t_0),
\xi^+)\\
&\ge \dis(\Sigma_p(S), \xi^-)+\dis(\Sigma_p(S), \xi^+)\\
&=\ell^-+\ell^+\\
&\ge \frac{\pi}{2}+\frac{\pi}{2}\\
&=\pi
\end{split}
\end{equation}
Note that $\curv(\Sigma_p(X^4))\ge 1$ implies
$\diam(\Sigma_p(X^4))\le \pi$, hence the inequalities in
\eqref{eq:1} are equalities, in particular $\ell^-=\ell^+=\pi/2$ and
$\Sigma_p(X^4)$ is the spherical suspension over $\Sigma_p(S)$, i.e.
$T_p(X^4)$ splits isometrically as $T_p(S)\times\RR^1$. Hence we
have a normal line bundle over $S$. By passing to the double cover,
we can assume that this line bundle is trivial, therefore $S$
separates $X^4$ into two parts and $X^4$ has two ends. Now it's easy
to see $X^4$ admits a line, by splitting theorem, $X^4$ is isometric
to $S\times \RR^1$.
\end{proof}

\subsection{Soul of the types $(s=0 \text { or } 2, m=3)$, $X^4\cong \RR^4$ or
$\RR^2\hookrightarrow X^4 \to \Sigma^2$ where $\Sigma^2\cong
\SS^2/\Gamma, \TT^2/\Gamma$}\label{sub:4D:2}
\begin{proof}
In this case, $\Omega(0)\cong D^3$ or $I$-bundle over $S\cong \SS^2,
\RP^2, \TT^2$ or $\KK^2$ by \autoref{thm:3D:3DSoul}. By the proof of
\autoref{sub:4D:1}, the interior of $\Omega(0)$ admits a normal line
bundle. Thus, we only have to show that $B_{X^4}(\Omega(0),
\varepsilon)\cong B_{X^4}(\Omega^{-\delta}(0), \varepsilon)$ for
$0<\varepsilon\ll \delta$, where $\Omega^{-\delta}(0):= \{x\in
\Omega| \dis(x, \partial \Omega(0))\ge \delta\}$.

By \autoref{lem:4D:nonCritcal}, $r_{\partial\Omega(0)}(x)=
\dis(\partial\Omega(0), x)$ has no critical point for $x\in
X^4\setminus\Omega(0)$. When restrict to $\Omega(0)$,
$r_{\partial\Omega(0)}$ is concave, hence it has no critical value
in $(0, a)$ for $a$ small enough. Combine these two facts, we know
there exits $\delta>0$ such that $r_{\partial\Omega(0)}$ has no
critical point in $B_{X^4}(\partial\Omega(0), 100\delta)$.

By the lower semi continuity of the norm of gradient of
$\lambda$-concave function $\nabla r_{\partial\Omega(0)}$  (cf.
\cite{Petr2007} Corollary 1.3.5), $r_{\partial\Omega(0)}$ has no
critical point in $B_{X^4} (\Omega(0) \setminus
\Omega^{-2\delta}(0), \varepsilon) - \partial\Omega(0)$ for
$\varepsilon\ll \delta$ small enough. Therefore
$r_{\partial\Omega(0)}$ is regular in
$B_{X^4}(\Omega^{-\varepsilon}(0) \setminus \Omega^{-2\delta}(0),
\varepsilon)$, by Fibration Theorem we have
\begin{equation}\label{eq:chap3:homeo1}
B_{X^4}(\partial\Omega(0), \varepsilon)\cong
B_{X^4}(\Omega(0)\setminus\Omega^{-2\delta}(0), \varepsilon)
\end{equation}
for $\varepsilon\ll \delta$, see \autoref{Pic:1}.
\begin{figure*}[ht]
% Requires \usepackage{graphicx}
\includegraphics[width=200pt]{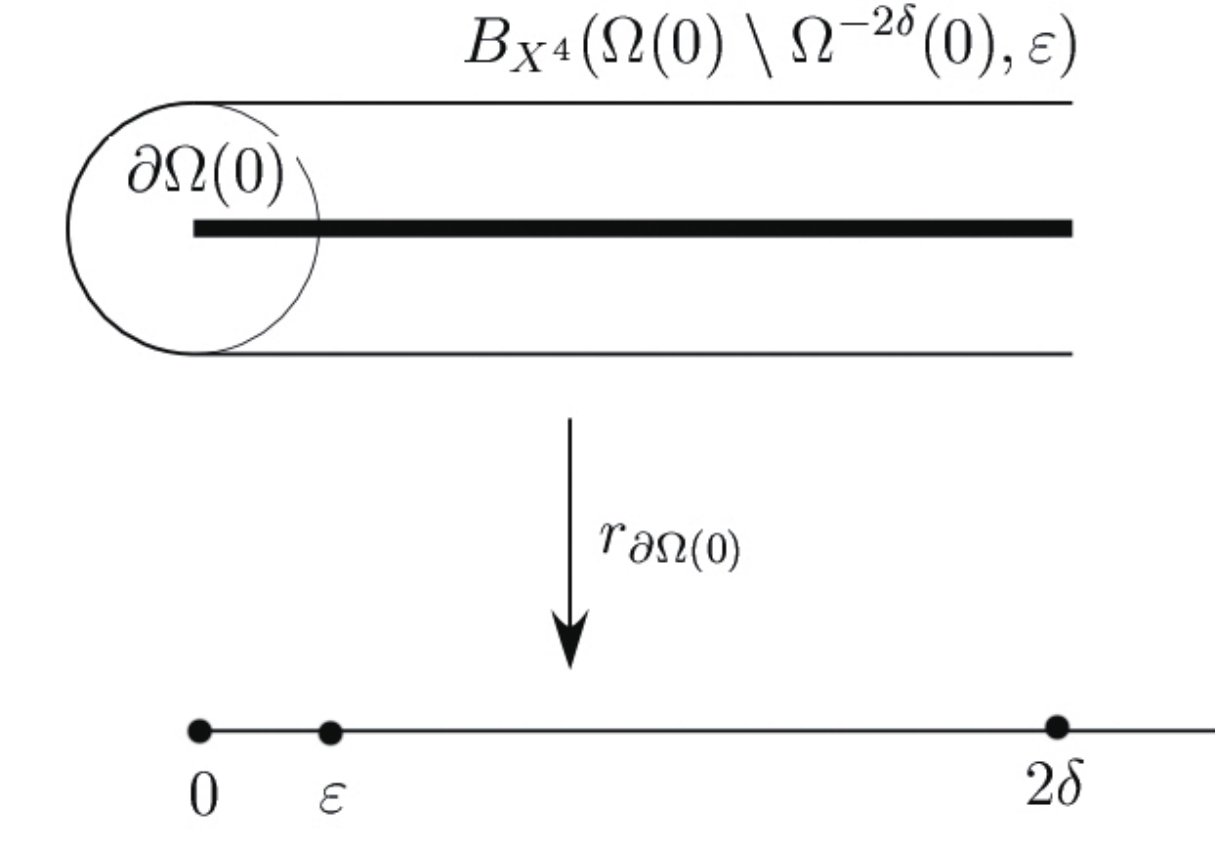}\\
\caption{Bundle structure around $\partial\Omega(0)$}\label{Pic:1}
\end{figure*}

By the proof of \autoref{sub:4D:3}, we know $B_{X^4}
(\partial\Omega(0), \varepsilon)$ is homeomorphic to a $D^2$ bundle
over $\partial \Omega(0)$
\begin{equation}\label{eq:chap3:bdyBdle}
D^2\hookrightarrow B_{X^4} (\partial\Omega(0), \varepsilon)\to
\partial \Omega(0)
\end{equation}
since $\partial\Omega^{-\varepsilon/2}(0)\cong \partial\Omega(0)$
and there is a normal line bundle over
$\partial\Omega^{-\varepsilon/2}(0)$, by passing to the double cover
one can assume this line bundle is trivial, let $\Gamma$ be the
$(\varepsilon/100)$-section of this line bundle. Clearly
$\Gamma\subset B_{X^4} (\partial\Omega(0), \varepsilon)$, which
implies that the bundle \eqref{eq:chap3:bdyBdle} admits a global
nowhere vanishing section, hence it is a trivial $D^2$ bundle. By
the homeomorphism \eqref{eq:chap3:homeo1}, we have a trivial $D^2$
bundle
\begin{equation}\label{eq:chap3:bdyBdle1}
B_{X^4}(\Omega(0)\setminus\Omega^{-2\delta}(0), \varepsilon)\cong
\partial\Omega(0)\times D^2
\end{equation}

Now consider the function $r_{\Omega^{-10\delta}}(x)$, which has no
critical point in $B_{X^4}(\Omega(0)\setminus\Omega^{-9\delta}(0),
\varepsilon) - \partial\Omega(0)$, hence we have a gradient flow of
$r_{\Omega^{-10\delta}(0)}$ on
\begin{equation}\label{eq:chap3:intBdle}
\begin{split}
B_{X^4}(\Omega^{-\delta}(0)\setminus\Omega^{-2\delta}(0),
\varepsilon) &\cong (-\varepsilon, +\varepsilon)\times \partial
\Omega(0)\times(\delta, 2\delta)\\
&\cong \partial \Omega(0)\times D^2
\end{split}
\end{equation}
for $\varepsilon\ll \delta$ small enough, where the homeomorphism
follows from the facts that $r_{\partial \Omega(0)}$ is concave in
the interior of $\Omega(0)$ and that there is a normal line bundle
over the interior of $\Omega(0)$, see \autoref{sub:4D:1}.

%Hence the gradient flow gives a homeomorphism:
%$$
%B_{X^4}(\Omega(0)\setminus\Omega^{-2\delta}(0), \varepsilon) \cong
%B_{X^4}(\Omega(0)\setminus\Omega^{-\delta}(0), \varepsilon)
%$$

Since the bundle \eqref{eq:chap3:bdyBdle1} and it's sub-bundle
\eqref{eq:chap3:intBdle} are both trivial bundles, one can extend
the gradient flow of $r_{\Omega^{-10\delta}(0)}$ on
$B_{X^4}(\Omega^{-\delta}(0) \setminus\Omega^{-2\delta}(0),
\varepsilon)$ to a gradient-like flow on
$B_{X^4}(\Omega(0)\setminus\Omega^{-2\delta}(0), \varepsilon)$,
which give the homeomorphism (see \autoref{Pic:2}):
$$
B_{X^4}(\Omega(0), \varepsilon) \cong B_{X^4}(\Omega^{-2\delta}(0),
\varepsilon)
$$

\begin{figure*}[ht]
% Requires \usepackage{graphicx}
\includegraphics[width=200pt]{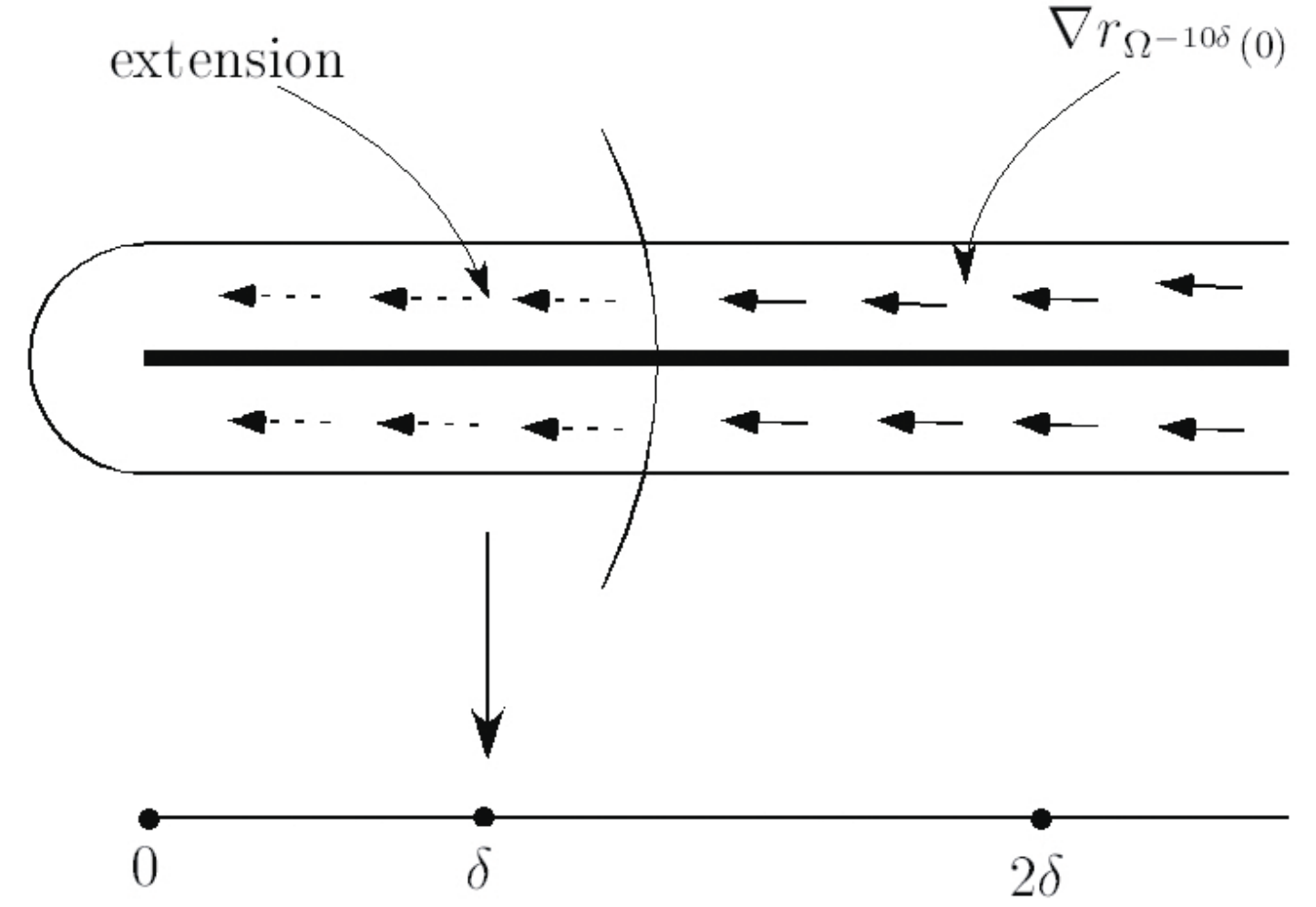}\\
\caption{Gradient-like flow on
$B_{X^4}(\partial\Omega(0),\varepsilon)$}\label{Pic:2}
\end{figure*}

It follows from the proof of \autoref{sub:4D:1} that $B_{X^4}
(\Omega^{-2\delta}(0), \varepsilon)$ is homeomorphic to the normal
line bundle over $\Omega^{-2\delta}(0)$, hence combine with
\autoref{thm:3D:3DSoul}, our main theorem holds in these two cases.

\end{proof}

\subsection{Soul of the type $(s=2, m=2)$,
$\RR^2\hookrightarrow X^4\to S$, $S=\SS^2/\Gamma,
\TT^2/\Gamma$}\label{sub:4D:3}

\begin{proof}
In this case, $S$ is $2$-dimensional surface with non-negatively
curvature. Thus $S$ is homeomorphic to $\SS^2, \TT^2, \RP^2$ or
$\KK^2$, by Splitting Theorem, if $S=\TT^2$  or $S=\KK^2$ the
universal cover $\widetilde{X^4}$ splits isometrically as $N^2\times
\RR^2$, thus the theorem follows from the fact that $N^2$ is
homeomorphic to $\RR^2$. Now we consider the cases where $S=\SS^2$
or $\RP^2$.

Let $\{p_i\}_{i=1}^N$ be the set of extremal points on $S$, where
the extremal points in Alexandrov surface is defined to be the
points satisfying $\diam(\Sigma_{p_i}(S))\le \pi/2$. Note that
$\Sigma_{p_i}(S)=\SS^1$ is a convex subset of
$\Sigma_{p_i}(X^4)=\SS^3$, hence for $\varepsilon>0$, we have
$B_{T_{p_i}(X^4)}(T_{p_i}(S), \varepsilon) = T_{p_i}(S)\times D^2$.
Then by Stability Theorem, we know there exits $\delta>0$ such that
$B_{X^4}(p_i, \delta)\cap B_{X^4}(S, \varepsilon)\cong D^4$ and for
$\varepsilon\ll \delta$, we have disc bundle

\begin{equation}\label{eq:chap3:Bundle1}
D^2\hookrightarrow B_{X^4}(p_i, \delta)\cap B_{X^4}(S,
\varepsilon)\xrightarrow{\pi_i} B_{S}(p_i, \delta).
\end{equation}
where $\pi_i$ is the bundle projection map. In particular
$\pi^{-1}_i(\partial B_S(p_i, \delta))= D^2\times S^1$

For $p\in S\setminus(\cup_{i=1}^N B_S(p_i, \delta/10))$, by our
assumption $\Sigma_p(S)> \pi/2$, thus there exits $\delta'>0$ and a
admissible map $F_p=(f_1, f_2): B_S(p, \delta')\to \RR^2$ which is
regular $B_S(p, \delta')$, by the lower semi-continuity of gradient
of semi-concave functions, we know $F_p$ is regular in $B_{X^4}(p,
\delta'')$, for some $\delta''>0$ satisfying $\delta''\le \delta \ll
\varepsilon$. Thus we have a fiber bundle:
\begin{equation}\label{eq:chap3:Bundle2}
N\hookrightarrow B_{X^4}(p, \delta'')\xrightarrow{\pi} D^2\cong
B_S(p, \delta'')
\end{equation}

Let $f(x)=\dis(x, S)$. Since $G_p=(f_1, f_2, f)$ is regular in the
domain $A_{X^4}(S, \varepsilon/100, \varepsilon)\cap B_{X^4}(p,
\delta'')$, for $\varepsilon\ll \delta$, where $A_X(S,
\varepsilon_1, \varepsilon_2)$ denoted the annular region, i.e. all
points have distance to $S$ between $\varepsilon_1$ and
$\varepsilon_2$. It follows from the Fibration Theorem that we have
a fiber bundle:
$$
\SS^1 \hookrightarrow A_{X^4}(S, \varepsilon/100, \varepsilon)\cap
B_{X^4}(p, \delta'')\to D^2\times I
$$
where $I$ is a open interval. Thus $\partial N=\SS^1$, by
generalized Margulis Lemma (cf. \cite{FY1992}), $N\cong D^2$.

Now we can glue the $D^2$ bundle together, (this part is similar to
Yamaguchi's construction in \cite{Yam2002} Page102). Let
$S(\frac{\delta}{2})= S-\cup_{i=1}^NB(p_i,\frac{\delta}{2})$. By
construction we have a $D^2$ bundle over $S(\frac{\delta}{2})$:
\begin{equation}\label{eq:chap3:bdloverreg}
D^2\hookrightarrow B_{X^4}(S(\frac{\delta}{2}),
\varepsilon)\xrightarrow{\pi} S(\frac{\delta}{2})
\end{equation}
Hence $\pi^{-1}(\partial B_S(p_i, \frac{\delta}{2}))=D^2\times S^1$.
Consider the gradient flow of $r_{p_i}(\cdot)=\dis(p_i, \cdot)$ on
$A_S(p_i, \delta/2, \delta)$, again by the lower semi-continuity of
$|\nabla r_{p_i}|$, for $\varepsilon\ll\delta$ small enough,
$r_{p_i}$ is regular in $B_{X^4}(A_S(p_i, \delta/2, \delta),
\varepsilon)$ hence provide a homeomorphism $\phi$ between $F_i:=
\pi^{-1}(\partial B_S(p_i, \frac{\delta}{2}))$ and
$G_i:=\pi^{-1}(\partial B_S(p_i, \delta))$. Clearly $\partial
B_S(p_i, \delta/2)$ isotopic to $\partial B_S(p_i, \delta)$ in
$B_S(p_i, 2\delta)$, hence when restricted to the boundaries
$\partial F_i=\SS^1\times \SS^1$ and $\partial G_i=\SS^1\times
\SS^1$ , $\phi$ is isotopic to the identity, therefore we can glue
the $D^2$ bundle together.
\end{proof}

\subsection{Soul of the types $(s=1, m=1,2,3)$, $\RR^3\hookrightarrow X^4\to
\SS^1$}\label{sub:4D:4}

\begin{proof}
If the soul is $\SS^1$, then the universal cover $\widetilde{X^4}$
admits a line by the totally convexity of soul, so $\widetilde{X^4}=
N^3\times \RR$, where $N^3$ is homeomorphic to $\RR^3$ by
\autoref{thm:3D:3DSoul}, thus \autoref{thm:intr:Main} holds in this
case.
\end{proof}

\subsection{Soul of the type $(s=0, m=2)$, $X^4\cong \RR^4$}\label{sub:4D:6}
\begin{proof}
Since the normal bundle over point soul $S$ is homeomorphic to
$\RR^4$, it's enough to show that $B_{X^4}(\Omega(0),
\varepsilon)\cong D^4$. It's clear that $\Omega(0)\cong D^2$ since
the soul is a point. By the proof of \autoref{sub:4D:3} and the fact
that $D^2$ is contractible, we see that all $D^2$ bundle over
$\Omega(0)\cong D^2$ is trivial, we have
$$
B_{X^4}(\Omega^{-100\varepsilon}, \varepsilon)\cong D^2\times
D^2\cong D^4
$$
We claim that $B_{X^4}(\partial \Omega(0), \varepsilon)\cong
\SS^1\times D^3$. Assume the claim first, by the proof of
\autoref{sub:4D:2}, the gradient of
$\dis_{X^4}(\Omega^{-10\delta}(0), \cdot)$ can be extend to
$B_{X^4}(\partial \Omega(0), \varepsilon)$, and will give the
homeomorphism from $B_{X^4}(\Omega^{-100\varepsilon},
\varepsilon)\cong D^4$ to $B_{X^4}(\Omega(0), \varepsilon)$.

{\it Proof of the Claim:} Let $\{p_i\}_{i=1}^N$ be the set of
extremal points in $\partial \Omega(0)$. By the Stability Theorem,
$B_{X^4}(p_i, \varepsilon)\cong D^4$. Let $\gamma_i$ be the boundary
curve connection $p_i$ to $p_{i+1}$ with the understanding that
$p_{N+1}=p_1$. One can assume $\gamma_i$ is short enough such that
$\dis_{p_i}$ has no critical point in $B_{X^4}(\gamma_i^{-\delta},
\varepsilon)$, where $\delta\gg\varepsilon$, $\gamma^{-\delta}_i$ is
the sub-curve of $\gamma_i$ defined by $\{x\in
\gamma_i|\dis(x,p_i)\ge \delta {\ and\ }\dis(x,p_{i+1})\ge
\delta\}$. Hence by Fibration Theorem it's a locally trivial fiber
bundle,
$$
N^3\hookrightarrow B_{X^4}(\gamma_i, \varepsilon)\to
\gamma_i^{-\delta}
$$
Since $B_{X^4}(p_i, \varepsilon)\cong D^4$, we have $S_{X^4}(p_i,
\delta)\cap B_{X^4}(\gamma^{-\delta}_i, \varepsilon) \cong D^3$,
hence $N^3\cong D^3$. This finishes the proof.
\end{proof}

\subsection{Soul of the type $(s=0, m=1)$, $X^4\cong \RR^4$}\label{sub:4D:7}
\begin{proof}
The proof is identically same as Subcase 3.1 of Theorem 2.21 in
\cite{CaoG2010}, we omit it here.
\end{proof}

\subsection{Soul of the type $(s=0, m=0)$, $X^4\cong \RR^4$}\label{sub:4D:8}
\begin{proof}
Let $p=S$ be a soul and $\{\alpha_n\}_{n=1}^{\infty}$ be a sequence
of number such that $\lim_{n\to \infty}\alpha_n=\infty$. It's clear
that $\lim_{n\to \infty}(\alpha_n X^4, p)= (T_p(X^4), O)$. Then it
follows by Stability Theorem that $B_{\alpha_nX^4}(p,\varepsilon)
\cong B_{T_pX^4}(O, \varepsilon)\cong \RR^4$. Since $\dis_p(x)$ has
no critical point in $X^4\setminus p$, we conclude that $X^4\cong
\RR^4$ by Perelman's Fibration theorem.
\end{proof}

This completes the proof of the Main Theorem.\\

{\bf Acknowledgement} The author is indebted to his advisor
Professor Jianguo Cao for his guidance and support. The author is
also grateful to Professor Karsten Grove for many useful
discussions.

\end{document}